\documentclass[psamsfonts]{article}

    \usepackage{latexsym}
    \usepackage{amsfonts}
    \usepackage{amsmath}
    \usepackage{amssymb}
    \usepackage{amscd}
    \usepackage{amstext}
    \usepackage{amsthm}

 \title{Furstenberg Transformations and Approximate Conjugacy
 \thanks{Supported by National Science Foundation
of USA. \,\,\,\,\,\,\,\,
 AMS 2000 Subject Classification Number: Primary 37A55 and 46L35.\,\,\,\,\,\,\,\,\,\,\,
Key Words: Simple $C^*$-algebras, Approximate Conjugacy\protect\\}}

 \author{
 Huaxin Lin\\
    Department of Mathematics\\
    University of Oregon\\
    Eugene, Oregon 97403-1222\\
    }
 \date{}
  
\begin{document}
\maketitle

    %Definitions of repeated phrases
    \newcommand{\CA}{$C^*$-algebra}
    \newcommand{\SCA}{$C^*$-subalgebra}
\newcommand{\aue}{approximate unitary equivalence}
    \newcommand{\ayue}{approximately unitarily equivalent}
    \newcommand{\mops}{mutually orthogonal projections}
    \newcommand{\hm}{homomorphism}
    \newcommand{\pisca}{purely infinite simple \CA}
    \newcommand{\andeqn}{\,\,\,\,\,\, {\rm and} \,\,\,\,\,\,}
    \newcommand{\QED}{\rule{1.5mm}{3mm}}
    \newcommand{\morp}{contractive completely
    positive linear map}
    \newcommand{\asmorp}{asymptotic morphism}
    % Definitions of symbols
    \newcommand{\arrow}{\rightarrow}
    \newcommand{\tdsum}{\widetilde{\oplus}}
    \newcommand{\pa}{\|}  % \| probably works better than \parallel as a
                             % norm sign.
    \newcommand{\ep}{\varepsilon}
    \newcommand{\id}{{\rm id}}
    \newcommand{\aueeps}[1]{\stackrel{#1}{\sim}}
    \newcommand{\aeps}[1]{\stackrel{#1}{\approx}}
    \newcommand{\dt}{\delta}
    \newcommand{\yu}{\fang}
    \newcommand{\ca}{{\cal C}_1}
\newcommand{\Ad}{{\rm ad}}
    \newcommand{\tr}{{\rm TR}}
    \newcommand{\N}{{\bf N}}
    \newcommand{\C}{{\bf C}}
    \newcommand{\Aut}{{\rm Aut}}
    \newcommand{\Tand}{\,\,\,\text{and}\,\,\,}
\newcommand{\T}{{\mathbb T}}
\newcommand{\R}{{\mathbb R}}
    \newtheorem{thm}{Theorem}[section]
    \newtheorem{Lem}[thm]{Lemma}
    \newtheorem{Prop}[thm]{Proposition}
    \newtheorem{Def}[thm]{Definition}
    \newtheorem{Cor}[thm]{Corollary}
    \newtheorem{Ex}[thm]{Example}
    \newtheorem{Pro}[thm]{Problem}
    \newtheorem{Remark}[thm]{Remark}
    \newtheorem{NN}[thm]{}
    \renewcommand{\theequation}{e\,\arabic{section}.\arabic{equation}}
    \newcommand{\rforal}{{\rm\,\,\,for\,\,\,all\,\,\,}}
\newcommand{\Z}{{\mathbb Z}}
\newcommand{\Q}{{\mathbb Q}}
\newcommand{\di}{{\rm dist}}
       \newcommand{\af}{\alpha}
       \newcommand{\bt}{\beta}
       \newcommand{\Om}{\Omega}
        \newcommand{\om}{\omega}
       \newcommand{\gm}{\gamma}
      \newcommand{\sm}{\sigma}
      \newcommand{\ud}{\underline}
       \newcommand{\beq}{\begin{eqnarray}}
        \newcommand{\eneq}{\end{eqnarray}}
\begin{abstract}

%In the short note we study Furstenberg transformations on $2$-torus.
Let $\af$ and
$\bt$ be  two Furstenberg transformations on $2$-torus associated
with irrational numbers $\theta_1,$ $\theta_2,$ integers $d_1, d_2$ and Lipschitz functions
$f_1$ and $f_2.$
We show that $\af$ and $\bt$ are approximately conjugate in a measure theoretical sense if (and only
if) $\overline{\theta_1\pm \theta_2}=0$ in $\R/\Z.$ Closely related to the classification of simple
amenable \CA s, we show that $\af$ and $\bt$ are approximately $K$-conjugate if (and only if)
$\overline{\theta_1\pm \theta_2}=0$ in $\R/\Z$ and $|d_1|=|d_2|.$ This
is also shown to be equivalent to that the associated crossed product \CA s are isomorphic.

\end{abstract}

\section{Introduction}

A celebrated result of Giordano, Putnam and Skau  (\cite{GPS}) states that two minimal Cantor systems are strong orbit
equivalent if (and only if) the associated crossed product \CA s are isomorphic which can also be
described by their
$K$-theory. Moreover, two such Cantor systems are topological orbit equivalent
if part of the $K$-theoretical information (namely the tracial range of the $K_0$-groups) of the
associated
\CA s are (unital) order isomorphic.
This note is an attempt to explore  its possible analogy in the case that the space is connected.

With the recent rapid development of  the classification of amenable simple \CA s of stable rank one
(\cite{EG},
\cite{EGL}, \cite{Lntaf}, \cite{DE}, \cite{Lnann} and \cite{Lnduke}, to name few),
it becomes possible to apply \CA\, theory to the study of minimal homeomorphisms on more general
spaces. Several versions of approximate conjugacy have been introduced and studied recently (see
\cite{LnKC},
\cite{LM1},
\cite{M},
\cite{LM2} and \cite{LM3}).  In \cite{LM2} and \cite{LM3},
minimal homeomorphisms on the product of Cantor set and the circle were studied. It
 is shown that if a set of $K$-theoretical information of the two minimal homeomorphisms on
the product of Cantor set and circle are the same then they are approximate
$K$-conjugate (and the converse also holds).

%GPS's theorem is out of the idea on the classical measure theoretical orbit equivalence result of
%Dyer.
One of the reasons that  Giordano, Putnam and Skau's work is so  successful  is that Cantor set is
totally disconnected. Perhaps the orbit equivalence for Cantor minimal systems may be viewed  as
%an
%equivalence relation 
something which lives between measure theory and topology. When the space $X$ is
connected, the situation is very different. For example, by a result of Sierpinski, for connected
spaces, two topological orbit equivalent minimal homeomorphisms are in fact  flip conjugate (see Prop.
5.5 of
\cite{LqP}). Therefore, for connected space, one should not expect that two minimal homeomorphisms are
 topological orbit equivalent if their associated crossed product \CA s have
order isomorphic  $K$-groups, or if  their associated crossed product \CA s are isomorphic.

We are interested in the following two questions:

{\bf Q1}:\,\,\, Let $X$ be a (connected) compact metric space and let $\af$ and $\bt $ be two  minimal
homeomorphisms on $X.$ Let $A_\af$ and $A_\bt$ be the associated crossed product \CA s. Suppose that
the tracial range of
$K_0(A_\af)$ and that of $K_0(A_\bt)$ are (unital) order isomorphic.
What can one say about the homeomorphisms $\af$ and $\bt?$

\vspace{0.1in}

{\bf Q2}:\,\,\, Let $X$ be a (connected) compact metric space and let $\af$ and $\bt $ be two  minimal
homeomorphisms $X.$ Let $A_\af$ and $A_\bt$ be the associated crossed product \CA s. Suppose that
$A_\af$ and $A_\bt$ are isomorphic  with additional $K$-theoretical information
of $(j_\af)_*$ and $(j_\bt)_*.$ What can one say about the homeomorphism $\af$ and
$\bt?$

(See the definition of $j_\af$ and $j_\bt$ in \ref{DAA} below. A clarification of this
question will be discussed in \ref{KC=})

\vspace{0.1in}

Giordano, Putnam and Skau's results  answered  both questions, namely, $\af$ and $\bt$ are topological orbit equivalent for the
{\bf Q1} and $\af$ and $\bt$ are strong orbit equivalent for {\bf Q2} under  the assumption that $X$ is
the Cantor set.

The results in \cite{LM1} \cite{M}, \cite{LM2} and \cite{LM3} suggest that answer to {\bf Q2}
should be that
$\af$ and $\bt$ are approximately $K$-conjugate, and for {\bf Q1}, $\af$ and $\bt$ are approximately
conjugate in a more measure theoretical sense. However, spaces studied in the above mentioned articles
are not connected. It is interesting to see answers to the question
{\bf Q1} and {\bf Q2} for any connected spaces (other than $\T$).

A classical example of minimal homeomorphisms on the 2-torus $\T^2$ was studied
by Furstenberg (\cite{F}). Let $\theta$ be an irrational number and $g: \T\to \T$ be
a continuous map. The Furstenberg transform $\af: \T^2\to \T^2$ is defined to
be
$\af(\xi,\zeta)=(\xi e^{i2\pi\theta},\zeta g(\xi))$ for $\xi\in \T$ and $\zeta\in \T$ with $g$ being
homotopically non-trivial.
One may write that $\af=(\xi e^{i2\pi\theta},\zeta \xi^d e^{i2\pi f(\xi)}),$ where $d$ is an integer
and
$f$ is a real continuous function in $C(\T).$
Denote this  $\af$ by $\Phi_{\theta, d, f}.$ It is proved that $\af$ is always minimal if
$d\not=0$ (\cite{F}).   It is also shown in
\cite{F} that if $g$ satisfies the Lipschitz condition then
$\af$ is also unique ergodic (\cite{F}).
This is one of the favorite examples and has been intensively studied (for example,
\cite{F}, \cite{J}, \cite{Rou}, \cite{Ko1}, \cite{Ko2} and \cite{OP}, to name a few).

It was conjectured by R. Ji (\cite{J}) that
$\Phi_{\theta, 1,0}$ is conjugate to $\Phi_{\theta,1,f}.$  By
considering quasi-discrete spectrum, counter-examples have been constructed by Rouhani (\cite{Rou})
that
$\Phi_{\theta, 1, 0}$ may not be flip conjugate to
$\Phi_{\theta, 1,f}.$

Let $\af=\Phi_{\theta_1, d_1,f}$ and $\bt=\Phi_{\theta_2, d_2, f_2}.$
 In this note, we first show
that, if $\overline{\theta_1\pm \theta_2}=0$ in $\R/\Z,$ then $\af$ and $\bt$ are approximately
conjugate in a measure theoretical sense (see \ref{Dconm} below). We also show that the converse is
true, i.e., if $\af$ and
$\bt$ are approximately conjugate in that sense, then $\overline{\theta_1\pm \theta_2}=0$ in $\R/\Z.$ Let
$A_\af=C(\T^2)\rtimes_{\af}\Z$ and $A_\bt=C(\T^2)\rtimes_{\bt}\Z$ be the associated crossed
products. At least in the case that $f_1$ and $f_2$ satisfy the Lipschitz condition, the ranges of
$K_0(A_\af)$ and $K_0(A_\bt)$ under the tracial map are the same, namely, $\Z+\Z(\theta_1)$
(we still assume that $\overline{\theta_1\pm \theta_2}=0$ in $\R/\Z$).  This result seems closer to that of the topological orbit
equivalence in the Cantor minimal systems.
%So this may be viewed the analog of topological
%orbit equivalence in the connected case.
%Indeed, in a subsequent paper, we will show that a version of
%this holds in great generality.

It has been recently proved that (\cite{LP}), in the unique ergodic cases, $A_\af$ and $A_\bt$ are
unital simple
\CA s with tracial rank zero. Therefore, by the classification of unital simple
amenable \CA s with tracial rank zero (see \cite{Lnduke}), $A_\af$ and $A_\bt$ are isomorphic as \CA s if and only if
$\overline{\theta_1\pm \theta_2}=0$ in $\R/\Z$ and $|d_1|=|d_2|.$
We will show that $\af$ and $\bt$ are approximately
$K$-conjugate if and only if $\overline{\theta_1\pm \theta_2}=0$ in $\R/\Z$ and $|d_1|=|d_2|.$
In the process, we also show that, when $f_1-f_2$ is  in a  dense subset  of real  part of $C(\T),$
$\Phi_{\theta, d, f_1}$ and $\Phi_{\theta, d, f_2}$ are actually conjugate (see \ref{C1} below).

The results in this note are very special.
% and the proof should not be regarded as something that may
%have any significant impact on minimal homeomorphisms on more general connected spaces.
However, it is our hope that this  special case will lead us to more interesting answers to the
questions {\bf Q1} and {\bf Q2} and serves as an invitation for further exploration.

{\bf Acknowledgements} This work was partially supported by a NSF grant.
The author would like to thank
N. C. Phillips for many useful conversation.

\section{The Main results}

\begin{Def}\label{Dconm}
{\rm Let $X$  and $Y$ be two compact metric spaces and let $\af: X\to X$ and
$\bt: Y\to Y$ be two
minimal homeomorphisms. Let
$T_\af$ and $T_\bt$ be sets of $\af$-invariant and $\bt$-invariant
normalized Borel measures, respectively. We say that $\af,\bt$ are approximately conjugate in the
sense (M1) if there exist two sequences of homeomorphisms $\sigma_n: X\to Y$ and
$ \gamma_n: Y\to X$ and affine homeomorphisms $\Lambda_1: T_\af\to T_\bt$ and
$\Lambda_2: T_\bt\to T_\af$
such that
\beq\label{DM1}
\lim_{n\to\infty}\sup_{\nu\in T_\bt}\nu(\{y\in Y: \di(\sigma_n\af\circ \sigma_n^{-1}(y), \bt(y))\ge
a\})=0,
\eneq
\beq\label{DM2}
\lim_{n\to\infty}\sup_{\mu\in T_\af}\mu(\{x\in X: \di(\gamma_n\circ\bt\circ \gamma_n^{-1}(x), \af(x))\ge
a\})=0
\eneq
for all $a>0,$
\beq\label{DM3}
\Lambda_1(\mu)(S)=\mu(\sigma_n^{-1}(S)),
\andeqn \Lambda_2^{-1}(\nu)(G)=\nu(\gamma_n^{-1}(G)) \,\,\, n=1,2,...
\eneq
for all Borel sets $S\subset Y,$ $G\subset X$ and for all $\mu\in T_\af,$ $\nu\in T_\bt.$ }
\end{Def}

A couple of remarks are in order.

(1)  Suppose that there exists a homeomorphism $\sigma: X\to Y$ such that
$\sigma\circ\af\circ\sigma^{-1}=\bt.$
Define $\Lambda:T_\af\to T_\bt$ by $\Lambda(\mu)(S)=\mu(\sigma^{-1}(S))$ for
all Borel sets $S\subset Y.$  It should be noted that
\beq\label{Lam}
\Lambda(\mu)(\bt(S))&=&\mu(\sigma^{-1}\circ \bt(S))\\ \nonumber
&=&\mu(\sigma^{-1}\circ\bt\circ\sigma\circ\sigma^{-1}(S))\\ \nonumber
&=&\mu(\af\circ \sigma^{-1}(S))=m(\sigma^{-1}(S))=\Lambda(\mu(S))
\eneq
for all Borel sets $S\subset Y.$ So $\Lambda(\mu)\in T_\bt.$ In particular, conjugate homeomorphisms
are approximately conjugate in the sense (M1).

 In general,
a sequence homeomorphisms  $\{\sigma_n\}$ does not preserve measures even though both (\ref{DM1}) and
(\ref{DM2}) hold. One could  have $\lim_{n\to\infty}\mu(\sigma_n(S))=0.$ Here we require, in the
definition, that  $\{\sigma_n\}$ has some consistent information on measure spaces. So in the
definition \ref{Dconm}, one should view that the conditions in  (\ref{DM3}) are  important part of
the definition.

(2) In the definition of \ref{Dconm}, put
$$
E_n(a)=\{y\in Y:\di(\sigma_n\af\circ \sigma_n^{-1}(y), \bt(y))\ge a\})
$$
for $a>0.$ Put
$$
S_n(a)=\{x\in X: \di(\sigma_n\circ \af(x),\bt\circ\sigma_n(x))\ge a\}.
$$
Then
$$
\sigma_n^{-1}(E_n(a))=S_n(a).
$$
By (\ref{DM1}) and  (\ref{DM4}),
$$
\lim_{n\to\infty}\sup_{\mu\in T_\af}\mu(S_n(a))=0.
$$

(3) It is an easy exercise that approximately conjugacy in the sense (M1) is an equivalence relation
among minimal homeomorphisms.

\begin{Def}\label{DF1}
{\rm Let $\theta$ be an irrational number and $g: \T\to \T$ be a continuous function. A Furstenberg
transform is a map $\af: \T^2\to \T^2$ defined by $\af((\xi, \zeta))=(\xi e^{i2\pi\theta}, \zeta
g(\xi))$ for $\xi\in \T$ and $\zeta\in \T$ with $g$ being of degree $d\not=0$ (the winding number is
$d$). There is a real function
$f\in C(\T)$ such that $g(\xi)=\xi^d\exp(i2\pi f(\xi))$ for $\xi\in \T.$  We will call $\af$ the
Furstenberg transform associated with irrational number $\theta,$  integer $d$ and function
$f.$ It will
also be denoted by
$\Phi_{\theta, d, f}.$

We are interested in the case that $(\T^2,\af)$ is uniquely ergodic. It is known that $\af$ is always
minimal (see \cite{F}). It is also shown in \cite{F} that $(\T^2,\af)$ is uniquely ergodic if $g$ has
Lipschitz property (or $f$ is Lipschitz).  The unique invariant measure is the product of the
normalized Lebesgue measure $m_2=m\times m.$ We fix the following metric on $\T^2:$
$$
\di((\xi,\zeta),(\xi',\zeta'))=\sqrt{|\xi{\overline{\xi'}}-1|^2+|\zeta'\overline{\zeta'}-1|^2},
$$
where $\xi, \xi', \zeta,\zeta'\in \T.$

 We will keep these notation throughout this note.

}
\end{Def}

\vspace{0.1in}

\begin{thm}\label{MT1}
Let $\af=\Phi_{\theta_1, d_1, f_1}, \bt=\Phi_{\theta_2, d_2, f_2}: \T^2\to \T^2$ be  unique ergodic
Furstenberg transforms.  Then the following are equivalent:

{\rm (1)} ${\overline{|\theta_1\pm \theta_2|}}=0$   in $\R/\Z;$

{\rm (2)} $\af$ and $\bt$ are approximately conjugate in the
sense (M1),
\end{thm}

In this case, one can make
\beq\label{DM201}
\lim_{n\to\infty}\sup\{ \di(\sigma_n\af\circ \sigma_n^{-1}(y), \bt(y)):y\in Y)\}
=0,
\eneq
\beq\label{DM202}
\lim_{n\to\infty}\sup\{ \di(\gamma_n\circ\bt\circ \gamma_n^{-1}(x), \af(x)) : x\in X\}=0,
\eneq
if one  does not insists that $\sigma_n$ and $\gamma_n$  to be continuous everywhere.
More precisely, we have the following

\begin{thm}\label{MT1+}

Let $\af=\Phi_{\theta_1, d_1, f_1}, \bt=\Phi_{\theta_2, d_2, f_2}: \T^2\to \T^2$ be unique ergodic
Furstenberg transforms.  Then each condition  {\rm (1)} or {\rm (2)} in \ref{MT1} is also equivalent
to the following:

{\rm (3)} There are sequences of Borel equivalences $\{\sigma_n\}$
%((\xi,\zeta))=(\sigma_n^{(1)}((\xi, \zeta)),
%\sigma_n^{(2)}(\xi,\zeta))$
and $\{\gamma_n\}$
%((\xi,\zeta))=(\gamma_n^{(1)}((\xi,\zeta)),\gamma_n^{(2)}(\xi, \zeta)))$
%and an affine homeomorphism $\Lambda: T_\af\to T_\bt$
such that
\beq\label{DM201+}
\lim_{n\to\infty}\sup\{ \di(\sigma_n\af\circ \sigma_n^{-1}(y), \bt(y)):y\in Y)\}
=0,
\eneq
\beq\label{DM202+}
\lim_{n\to\infty}\sup\{ \di(\gamma_n\circ\bt\circ \gamma_n^{-1}(x), \af(x)) : x\in X\}=0,
\eneq
 and, for each $n,$ there exists a closed subset $F_n\subset \T$ and $K_n\subset \T$ such that
$m(F_n)=0$ and $m(K_n)=0,$ and
$\sigma_n$ and $\gamma_n$ are continuous on $\T\times (\T\setminus F_n)$ and on $\T\times (\T\setminus K_n)$
respectively, moreover,
\beq\label{DM204}
m(\sigma_n(S))=m(S),\,\,\,\,\,\
m(\gamma_n(G))=m(G)\,\,\,n=1,2,...
\label{DM4}
\eneq
for all Borel sets $S\subset \T^2.$
\end{thm}

\vspace{0.1in}

\begin{Def}\label{DTR}
Let $A$ be a stably finite unital \CA\, and let $T(A)$ be the tracial state space. Denote by
$Tr$ the usual (non-normalized) trace on $M_k.$ Define
$\rho: K_0(A)\to Aff(T(A))$ by $\rho([p])=\tau(p)$ for projections in $M_k(A),$ where $\tau=t\otimes
Tr$ and
$t\in T(A).$
\end{Def}

\begin{Def}\label{DAA}
{\rm Let $X$ be a compact metric space and let $\af: X\to X$ be a minimal homeomorphism. Then the
transformation group \CA, the crossed product, $C(X)\rtimes_\af \Z$ will be denoted by $A_\af.$

We will use $j_\af: C(X)\to A_\af$ for the natural embedding.

For a unital \CA\, $A,$ we denote by ${\rm ad}u\, (a)=u^*au$ for all $a\in A.$

We fix a unitary
$u_\af$ so that ${\rm ad}\,u_\af j_\af(f)=j_\af(f\circ \af)$ for all $f\in C(X).$
It should be noted that there are other choices for such $u_\af.$  For example, if $z\in C(X)$ is a unitary
then  $w=u_\af j_\af(z)$ is another choice. In fact,
\beq\label{DAA1}
{\rm ad}\, w (j_\af(f))=j_\af(z)^*u_\af^*  (j_\af(f))u_\af j_\af(z)=j_\af(z^* (f\circ \af)
z)=j_\af(f\circ\af)
\eneq
for all $f\in C(X).$

}

\end{Def}

\begin{Remark}\label{RCK}
{\rm Let $X=\T^2$ and let $\af=\Phi_{\theta, d, f}$ be a Furstenberg transformation with Lipschitz
$f.$ Then $A_\af$ is a unital simple \CA\, with a unique tracial state.

 It is
computed (see Example 4.9 of \cite{Ph}) that $K_0(A_\af)\cong \Z\oplus \Z\oplus \Z$ with
$\rho(K_0(A_\af))=\Z\oplus
\Z(\theta)\subset \R$  and
$$
K_0(A_\af)_+=\{m_1+m_2+m_3\in \Z^3: m_1+m_3\theta>0 \,\,\,{\rm or}\,\,\,
m_1=m_2=m_3=0\},
$$
where the first two copies of $\Z$ is identified with the image
of $K_0(C(\T^2))$ under the embedding $(j_\af)_{*0},$
and $K_1(A_\af)\cong \Z\oplus \Z\oplus \Z/d\Z \oplus \Z,$
where $\Z/d\Z\oplus \Z$ is the image of $K_1(C(\T^2))$ under $(j_\af)_{*1}.$  Let $z_1, z_2: C(\T^2)\to
\T$ be the functions defined by $z_1((\xi, \zeta))=\xi$ and $z_2((\xi,\zeta))=\zeta$ ($\xi, \zeta\in
\T$). Then $(j_\af)_{*1}([z_1])$ is the standard generator of $\Z/d\Z$ and
$(j_\af)_{*1}([z_2])$ is the standard generator of $\Z.$

Let $X$ be a compact metric space, and let $\af, \bt: X\to X$ be two minimal homeomorphisms. It is
certainly desirable to have two sequences of homeomorphisms $\{\sigma_n\}$ and $\{\gamma_n\}$ on $X$
so that
\beq\label{DKKC1}
\lim_{n\to\infty}\sup\{\di(\sigma_n\circ \af\circ \sigma_n^{-1}(x),\bt(x)): x\in X\}=0\andeqn
\eneq
\beq\label{DKKC2}
\lim_{n\to\infty}\sup\{\di(\gamma_n\circ \bt\circ \gamma_n^{-1}(x), \af(x)): x\in X\}=0
\eneq

However, when a sequence of maps (such as  $\{\sigma_n\}$ and $\{\gamma_n\}$) involved,
one also expects that the maps have something in common.
At least in the Cantor set case, without any consistency on the conjugating maps
$\{\sigma_n\}$ and $\{\gamma_n\},$ (\ref{DKKC1}) and (\ref{DKKC2}) are not so interesting
(see \cite{LM1}).
Even though
one should not expect these maps will converge in any meaningful way, one hopes that
some information about these maps are independent of $n.$
For example, one would like to require that both sequences preserve the measures as in
\ref{Dconm} as well as in (3) of Theorem \ref{MT1+}. To be more topologically interesting,
one may require that both sequences preserve some topological data. For example, in
\cite{LM1}, approximate $K$-conjugacy requires that both sequences preserve $K$-theory (in the
crossed products). We use the following definition in this note:

\begin{Def}\label{DKKC}
{\rm Let $X$ be a compact metric space, and let $\af, \bt: X\to X$ be two minimal homeomorphisms. Two
homeomorphisms $\af$ and $\bt$ are said to be {\it approximately $K$-conjugate} if there exist two
sequences of homeomorphisms $\{\sigma_n\}$ and $\{\gamma_n\}$ on $X$ such that (\ref{DKKC1}) and
(\ref{DKKC2})  hold, and there exists an isomorphism $\phi: A_\af\to A_\bt$ and sequences of unitaries
$u_n\in A_\bt$ and $v_n\in A_\af$ such that
\beq\label{DKKC3}
\lim_{n\to\infty}\|{\rm ad}\, u_n\circ \phi(j_\af(f))-j_\bt(f\circ \sigma_n^{-1})\|=0
\eneq
for all $f\in C(X)$ and
\beq\label{DKKC4}
\lim_{n\to\infty}\|{\rm ad}\, u_n\circ \phi(u_\af) -u_\bt z_n\|=0,
\eneq
where $z_n\in U_0(A_\bt)$ such that
\beq\label{DKKC5}
\lim_{n\to\infty}\|z_nj_\bt (f)-j_\bt(f)z_n\|=0
\eneq
for all $f\in C(X),$ and

\beq\label{DKKC6}
\lim_{n\to\infty}\|{\rm ad}\, v_n\circ \phi(j_\bt(f))-j_\af(f\circ \gamma_n^{-1})\|=0
\eneq
for all $f\in C(X)$ and
\beq\label{DKKC7}
\lim_{n\to\infty}\|{\rm ad}\, v_n\circ \phi(u_\bt) -u_\af y_n\|=0,
\eneq
where $y_n\in U_0(A_\af)$ such that
\beq\label{DKKC8}
\lim_{n\to\infty}\|y_nj_\bt (f)-j_\af(f)y_n\|=0
\eneq
for all $f\in C(X).$ }
\end{Def}

\begin{Remark}\label{Rdkc}
{\rm

(1) It should be noted that (\ref{DKKC1} and (\ref{DKKC3}) imply both (\ref{DKKC4}) and (\ref{DKKC5}).
In fact one can take
% imply that
%\beq\label{DKKC9}
%\lim_{n\to\infty}\|{\rm ad}\, (u_n^*u_\af^*u_n) (j_\bt(f))-j_\bt (f\circ \bt)\|=0.
%\eneq
%This was used in \cite{LM1}, \cite{LM2} and \cite{lM3} instead of (\ref{DKKC4}) and (\ref{DKKC6}. Let
$z_n=u_\bt^*(u_n^*u_\af u_n).$
%Then (\ref{DKKC9}) implies that both (\ref{DKKC4}) and (\ref{DKKC5}).
At least from the discussion
that leads to (\ref{DAA1}), one sees that the unitaries $z_n$ and $y_n$ can not be omitted.

(2) The existence of the unitaries $\{u_n\}$ and $\{v_n\}$ implies that $\{\sigma_n\}$ and
$\{\gamma_n\}$ preserve the invariant measures.
Note also if $p\in M_k(C(X))$ is a
projection, then $[j_\bt (p\circ \sigma_n^{-1})]=[j_\bt(p)]$ in $K_0(A_\bt)$ for
all large $n.$  In fact, $\{\sigma_n\}$ and
$\{\gamma_n\}$ preserve
the ordered $K$-theory (independent of $n$) and beyond.

 }
\end{Remark}

Condition (1) in Theorem \ref{MT1} implies that $\rho(A_\af)$ and $\rho(A_\bt)$ are (unital) order
isomorphic. In the case that $X$ is the Cantor set, by a Giordano, Putnam and Skau's theorem, the
condition that
$\rho(A_\af)$ and $\rho(A_\bt)$ are unital order isomorphic is equivalent to that
$\af$ and $\bt$ are  topological orbit equivalent. Similar conclusion
is not possible for connected space, as mentioned earlier  that two topological orbit equivalent
minimal homeomorphisms on connected spaces are flip  conjugate. As we stated in the introduction,
topological orbit equivalence for minimal Cantor  systems seems to be something between measure theory
and topology. The conclusion of Theorem \ref{MT1} also has both topological and measure theoretical
favor. It should be noted that one can not make maps
$\{\sigma_n\}$ and
$\{\gamma_n\}$ in (3) of Theorem \ref{MT1+} to be homeomorphisms in general
as the next theorem states.
 It should also be mentioned that the approximate conjugacy in the sense of
(5) in Theorem
\ref{MT2} is rather  weak relation for minimal Cantor system (see
\cite{LM1} ).  However, this notion plays completely different role for
homeomorphisms on the connected spaces.

It is recent proved (\cite{LP}) that for  a unique ergodic Furstenberg transformation $\af$ the
associated crossed product $A_\af$ has tracial rank zero. So classification theorem (see
\cite{Lnann} and \cite{Lnduke}) can be
applied. In particular, if $\af=\Phi_{\theta_1, d_1, f_1}$ and $\bt=\Phi_{\theta_2, d_2, f_2},$ then
$A_\af\cong A_\bt$ is and only if $\overline{\theta_1\pm\theta_2}=0$ in $\R/\Z$ and $|d_1|=|d_2|.$

}

\end{Remark}

\begin{thm}\label{MT2}
%Let $\theta_1$ and $\theta_2$ be two irrational and $f_1, f_2: \T\to \R$ be continuous functions
%satisfying the Lipschitz condition.  Suppose that
%$\af=\Phi_{\theta, d_1,f_1}$ and $\bt=\Phi_{\theta, d_2, f_2}$
%are two Furtsenberg transforms.
Let $\af, \bt: \T^2\to \T^2$ be two unique ergodic Furstenberg transforms associated with irrational
numbers $\theta_1,\,\theta_2$ and integers $d_1, d_2\in \Z\setminus \{0\},$ respectively.

Then the following are equivalent.

{\rm (1)} $A_\af\cong A_\bt;$

{\rm (2)}  $\overline{|\theta_1\pm \theta_2|}=0$ in $\R/\Z$ and
$|d_1|=|d_2|;$

{\rm (3)}  $\af$ and $\bt$ are approximately $K$-conjugate;

{\rm (4)}  $(K_0(A_\af), K_0(A_\af)_0, [1_{A_\af}], K_1(A_\af))\cong (K_0(A_\bt), K_0(A_\bt)_+,
[1_{A_\bt}], K_1(A_\bt)).$

{\rm (5)} There exist two sequences of homeomorphisms $\{\sigma_n\}$ and $\{\gamma_n\}$ on
$\T^2$ such that

\beq\label{T2S1}
\lim_{n\to\infty}\sup\{\di(\sigma_n\circ \af\circ \sigma_n^{-1}(x),\bt(x)):x\in \T^2\}=0,
\eneq
\beq\label{T2S2}
\lim_{n\to\infty}\sup\{\di(\gamma_n\circ \bt\circ \gamma_n^{-1}(x), \af(x)):x\in
\T^2\}=0,
\eneq
\beq\label{T2S3}
m_2(\sigma_n(S))=m_2(S)\andeqn m_2(\gamma_n(S))=m_2(S)
\eneq
for all Borel subsets
$S\subset \T^2.$

\end{thm}

\begin{Remark}\label{KC=}
{\rm  It should be pointed out, in general, (1) or (4)
do not imply that $\af$ and $\bt$ are approximately K-conjugate.  See, for example, Example 9.2
of \cite{LM2}.
Let $X$ be a connected compact CW -complex
and let $\af, \bt: X\to X$ be two minimal homeomorphisms such that
$(X,\af)$ and $(X, \bt)$ are unique ergodic.
The right
condition,  in this case, should be, at least in the case that $\rho(K_0(A_\af))$ is dense in $K_0(A_\af),$ the following:
There is an order isomorphism
$$\kappa: (K_0(A_\af), K_0(A_\af)_+, [1_{A_\af}], K_1(A_\af))
\to (K_0(A_\bt)), K_0(A_\bt)_+, [1_{A_\bt}], K_1(A_\bt)),$$
there exists a sequence of  isomorphisms $\lambda_n:C( X)\to C(X)$ such that
$$
[\kappa\circ j_\af\circ  \lambda_n)]=[j_\bt] \,\,\,{\rm in}\,\,\, KL(C(X), A_\bt)
$$
and
$$
\lim_{n\to\infty}\tau\circ j_\af(f\circ \lambda_n)=\tau\circ j_\af(f),
$$
where $\tau$ is the unique tracial state of $A_\af$
(for more general case, see \cite{LnKC}). }
\end{Remark}

\section{Measure-theoretically approximate conjugacies}

The following is taken from p. 299 of \cite{LqP}  (see also Theorem 2.3 of \cite{LP}). 

\begin{thm}\label{Rok1}
Let $X$ be an infinite compact metric space, and let $h \colon X \to X$ be a minimal homeomorphism.
Let
$Y
\subset X$ be closed and have nonempty interior. For $y \in Y$ set $r (y) = \min \{ n \geq 1 \colon
h^n (y)
\in Y
\}.$ Then $\sup_{y \in Y} r (y) < \infty.$ Let $n(0)<n (1) < n (2) < \cdots < n (l)$ be the distinct values
in the range of $r,$ and for $0 \leq k \leq l$ set
\[
Y_k = {\overline{\{ y \in Y \colon r (y) = n (k) \} }}
\andeqn
Y_k^{\circ} = {\rm int}(\{ y \in Y \colon r (y) = n (k) \}).
\]
Then:
\begin{itemize}
\item[{\rm (1)}]
The sets $h^j (Y_k^{\circ}),$ for $0 \leq k \leq l$ and $1 \leq j \leq n (k),$ are disjoint.
\item[{\rm (2)}]
$\bigcup_{k = 1}^l \bigcup_{j = 1}^{n (k)} h^j (Y_k) = X.$
\item[{\rm (3)}]
$\bigcup_{k = 0}^l h^{n (k)} (Y_k) = Y.$
\end{itemize}
\end{thm}

\begin{Lem}\label{Rok2}
Let  $\af: \T\to \T$ be  a minimal homeomorphism. Let $n>1$ be an integer. Then there are finitely
many pairwise disjoint open arks
$\{J_i=(c_i,d_i): i=1,2,...,k\}$ of $\T$ such that

{\rm (1)} $\af^j(J_i)$ are pairwise disjoint for $0\le j\le h(i)-1$ and $i=1,2,...,k;$

{\rm (2)} $n\le h(i),$ $1\le i\le k;$

{\rm (3)} $\T\setminus \cup_{i=1}^k \cup_{j=0}^{h(i)-1}\af^j(J_i))$ is a set of finite many points.

\end{Lem}

\begin{proof}

Fix $x\in \T$ and fix an integer $n>1.$ Since $\af$ is minimal, there is a closed ark $Y$ containing
$x$ such that $\af^j(J)$ are pairwise disjoint for $0\le j\le n.$ Set
$$
r(y)=\min\{m\ge 1: \af^m(y)\in Y\}.
$$

Applying \ref{Rok1}, we obtain $n(0)<n(1)<\cdots n(l)$ and
$Y_0, Y_1,..., Y_l$ and
$Y_0^o, Y_1^o,...,Y_l^o$ as in the Lemma \ref{Rok1}. Note that $n(0)\ge n.$

Let $X_k=\{y\in Y: r(y)=n(k)\},$ $k=0,1,2...,l$ and let $\Om={\rm int} Y.$
Note that $X_0=Y_0.$  Let
$$
V_1=\af^{n(0)}(\Om)\cap \Om.
$$
Since $\Om$ is a non-empty   ark, so is $V_1.$ Set $S_1=\af^{-n(0)}(V_1).$ Then $S_1$ is an open
sub-ark of
$Y.$
Note that $S_1=Y_0^o.$

Let
$$
V_2=\af^{n(1)}(\Om)\cap \Om.
$$
Then $V_2$ is a non-empty open ark. Set  $S_2=\af^{-n(1)}(V_2).$ Then $Y_1^o=S_2\setminus X_0.$ This
shows that $X_1$ is a union of finitely many arks (possibly neither open nor closed).

Let $V_3=\af^{n(2)}(\Om)\cap \Om.$ Then $V_3$ is a non-empty open ark. Set $S_3=\af^{-n(2)}(V_3).$ Then
$Y_2^o=S_3\setminus (X_1\cup X_2).$
So $X_2$ is a union of finitely many arks.

By induction, we conclude that  all $X_k$ are union of finitely many arks. It follows that
$Y_k^o$ is a union of finitely many open arks and $Y_k$ is the closure of $Y_k^o.$
Note that
$\T\setminus \cup_{k=0}^l\cup_{j=0}^{n(j)}Y_k^o$ contains only finitely many points.
The lemma then follows.
\end{proof}

A similar lemma like below for the Cantor set appeared in the proof of 4.4 of \cite{M}. The following
is a version of that for the circle. Note the function $\omega$ can not be made to be continuous
on the whole space  as in
the Cantor set case.

\begin{Lem}\label{Lconj}
Let $\{I_1,I_2,...,I_k\}$ be  finitely many disjoint open arks of
$\T$ and let $\af:\T\to \T$ be a homeomorphism  such that $\af^j(I_i)$ are
pairwise disjoint for $0\le j\le h(i)-1$ and $1\le i\le k.$  Let $F, G:\T\to \T$ be  continuous maps.
Then, for any $\ep>0,$ and for each $i,$ there are $I_i^{(s)}\subset I_i$ $s=1,2,...,m(i)$ disjoint open  arks
and there is a continuous maps $\om: S\to \T,$ where $S=\cup_{j=0}^{h(i)-1}
\cup_{i=1}^k\cup_{s=1}^{m(i)}\af^j(I_i^{(s)}),$ such that

{\rm (1)} $\om(x)=0$ if $x\in \cup_{s,i} I_j^{(s)}$

{\rm (2)}
\beq\label{Lc1}
|[F(\af^j(x))+\om(\af^j(x))]-[\om(\af^{j+1}(x))+G(\af^j(x))]|<{1\over{h(i)}}
\eneq
for all $x\in I_i^{(s)},$ $s=1,2,...,m(i),$ $j=0,1,...,h(i)-2$  {\rm (} we identify $\T$ with $\R/\Z$ {\rm )} and

{\rm (3)} in $\R/\Z,$
\beq\label{Lc2}
|[F(\af^{h(i)-1}(x))+\om(\af^{h(i)-1}(x))]-[G(\af^{h(i)-1}(x))]|<{1\over{h(i)}}
\eneq
for all $x\in I_i^{(s)},$ $s=1,2,...,m(i),$ $i=1,2,...,k,$ and

{\rm (4)}
$I_i\setminus \cup_{s=1}^{m(i)}I_i^{(s)}$ contains only finitely many points.
%$m(\cup_{j=0}^{h(i)-1}\cup_{i=1}^k\af^j(I_i)\setminus (
%\cup_{j=0}^{h(i)-1}\cup_{i=1}^k\cup_{s=1}^{m(i)}\af^j(I_i^{(s)})))<\min\{\ep, d\},$
%where
%$$
%d=(1/2)\min\{m(I_i^{(s)}): 1\le s\le m(i),\,1\le i\le k\}.
%$$

Moreover,  on the closure of each $I_i^{(s)},$ $\omega$ can be extended
to be a continuous function and {\rm (1)}, {\rm (2)} and {\rm (3)} remain
true for $x$ in the left-closed arks of $I_i^{(s)} $ if inequalities
are replaced by  $``\le".$
\end{Lem}

\begin{proof}
We identify $\T$ with $\R/\Z.$ Define
\beq\label{LC1}
\kappa(x)=\sum_{i=0}^{h(i)-1} (F(\af^j(x))-G(\af^j(x)))
\eneq
for $x\in I_i,$ $1\le i\le k.$

Let $K=\sum_{i=1}^k h(i).$ One can break $I_i$ into a disjoint union of finitely many arks so that the
image of $\kappa$ on each sub-ark is a subset of a proper closed subset of $\T.$ Thus, there are, for
each
$i,$ pairwise disjoint open sub-arks $I_i^{(s)}\subset I_i,$
$s=1,2,..,m(i),$ such that $\kappa|_{{\overline{I_j^{(s)}}}}$ has image contained in a proper subset
of $\T$ such that
$$
I_i\setminus \cup_{s=1}^m(s)I_i^{(s)}
$$
contains only finitely many points, $i=1,2,...,k.$
So  (4) now holds.

Define ${\tilde \kappa}(x)$ as follows
\beq\label{LC2}
{\tilde \kappa}(x)=\kappa(x)+\Z \andeqn -1\le  {\tilde \kappa}(x)\le 1
\eneq
for $x\in \overline{I_i^{(s)}},$ $1\le s\le m(i),$ $1\le i\le k.$

Put $S=\cup_{i=1}^k\cup_{s=1}^{m(i)}I_i^{(s)}.$
Define $\eta: S\to \T$ as follows:
\beq\label{LC3}
\eta(\af^j(x))=-{j\over{h(i)}}{\tilde \kappa}(x)+\Z
\eneq
for $x\in \overline{I_i^{(s)}},$ $j=0,1,...,h(i)-1,$ $1\le i\le k.$
Since the image of $\kappa$ on each $\overline{I_i^{(s)}}$ is a proper subset of $\T,$
we see that $\eta$ is continuous.

Put $\Om_0=\sum_{i=1}^k\cup_{s=1}^{m(i)}I_i^{(s)}.$
Define
$\om: S\to \T$ as follows:
\beq\label{LC4}
\om(x)=0
\eneq
if $x\in \Om_0,$ and
\beq\label{LC5}
\om(\af^j(x))=\eta(\af^j(x))+\sum_{l=0}^{j-1}[F(\af^l(x))-G(\af^l(x))]
\eneq
for $x\in I_i^{(s)},$ $j=1,2,...,h(i)-1,$ $1\le s\le m(i),$  $1\le i\le k.$

If $x\in I_i^{(s)},$
\begin{align}\label{LC6}
&|[F(x)+\om(x)]-[\om(\af(x))+G(x)]|\\
&\le |F(x)-[{1\over{h(i)}}({\tilde \kappa}(x)+(F(x)-G(x))+G(x)]|<{1\over{h(i)}}
\end{align}
for $1\le s\le m(i)$ and $1\le i\le k.$

If $x\in \af^{j}(I_i^{(s)})$ for $j=1,2,...,h(i)-2,$
\begin{align}\label{LC7}
&|[F(\af^{j}(x))+\om(x)]-[\om(\af^{j+1}(x))+G(\af^{j}(x))]|\\
&\le |\eta(\af^j(x))-\eta(\af^{k+1}(x))|<{1\over{h(i)}}
\end{align}
for $1\le i\le k.$

We verify that, in $\R/\Z,$
\begin{align}\label{LC8}
&|[F(\af^{h(i)-1}(x))+\om(\af^{h(i)-1}(x))]-[G(\af^{h(i)-1}(x))]|\\
&=|\kappa(x)-{h(i)-1\over{h(i)}}{\tilde \kappa}(x)|<{1\over{h(i)}}
\end{align}
for all $x\in I_i,$ $1\le i\le k.$

We also note that the last statement follows easily since $F$ and $G$ are
continuous functions.
\end{proof}

\begin{Lem}\label{MINV}
Let $\sigma((\xi, \zeta))=(\sigma_1(\xi, \zeta), \sigma_2(\xi,\zeta))$
be a Borel equivalence from $\T^2$ to $\T^2$ such that
\beq\label{MINV1}
m_2(\sigma(S))=m_2(\sigma(S))
\eneq
for all Borel set $S.$

{\rm (1)} Then, for any Borel set $S_1, S_2\subset \T,$
\beq\label{MINV2}
m(\sigma_1((S_1,\zeta))=m(S_1)\andeqn m(\sigma_2((\xi, S_2))=m(S_2)
\eneq
for almost all $\zeta\in \T$ and almost all $\xi\in \T,$

{\rm (2)}  If there exists a closed subset $F\subset \T$ with $m(F)=0$ such that  $\sigma_1(-, \zeta)$
is continuous on $\T$ for all $\zeta\in \T,$ then
for each $\zeta\in \T,$
\beq\label{MINV3}
\sigma_1(\xi, \zeta)=\xi g_1(\zeta)
\eneq
for all $\xi\in \T$ or
\beq\label{MINV3+}
\sigma_1(\xi, \zeta)={\overline{\xi}}g_1(\zeta)
\eneq
for all $\xi\in T.$

{\rm (3)}
If $\sigma$ is continuous, then there are continuous maps $g_1: \T\to \T$ such that
\beq\label{MINV4}
\sigma_1(\xi, \zeta)=\xi g_1(\zeta)
\eneq
for all $\xi,\zeta\in \T$ or
\beq\label{MINV5}
\sigma_1(\xi,\zeta)={\overline{\xi}}g_1(\zeta)
\eneq
for all $\xi, \zeta\in \T.$

\end{Lem}
%If there exists a closed subset $F\subset \T^2$ such that
%$m_2(F)=0$ and $\sigma$ is continuous on $\T^2\setminus F,$
%then, there exists a measurable subset $E\subset \T$ with $m(\T\setminus E)=0$
%satisfying the following:  for each $\zeta\in  E,$ there exists pairwise disjoint
%open arks $O_i(\zeta)$ ($i=1,2,....$) such that, on each $O_i(\zeta),$
%\beq\label{MINV3}
%\sigma_1((\xi, \zeta))=\xi g_1(\xi, \zeta)\,\,\,{\rm or}\,\,\, \sigma_1((\xi, \zeta))=\overline{\xi}g_1(\xi, \zeta),
%\eneq
%where $g_1$ is a Borel function from $\T^2$ to $\T$ and is constant on each $O_n(\zeta).$

\begin{proof}

{\rm (1)} Let $S_1\subset \T$ be a Borel subset. Suppose that there exists a measurable subset $E_1\subset \T$
with positive measure such that
\beq\label{PMinv1}
\int_{S_1}\sigma_1(\xi, \zeta)d\xi\not=m(S_1)
\eneq
for all $\zeta\in E_1.$
Put
$$
E_1^+=\{\zeta\in E_1: \int_{S_1}\sigma_1(\xi, \zeta)d\xi >m(S_1)\}
\andeqn\\
E_1^-=\{\zeta\in E_1: \int_{S_1}\sigma_1(\xi, \zeta)d\xi <m(S_1)\}.
$$
If $m(E_1^+)>0,$ then
\beq\label{PMinv2}
m_2(S_1\times E_1^+)=\int_{S_1\times E_1^+} \sigma_1(\xi, \zeta)d\xi d\zeta
>\int_{E_1^+} m(S_1) d\zeta=m_2(S_1\times E_1^+).
\eneq
If $m(E_1^-)>0,$ then
\beq\label{PMinv3}
m_2(S_1\times E_1^-)=\int_{S_1\times E_1^-} \sigma_1(\xi, \zeta)d\xi d\zeta<\int_{E_1^-} m(S_1) d\zeta
=m_2(S_1\times E_1^-).
\eneq
Neither could be true. The proof for the variable $\zeta$ is the same.

{\rm (2)}
Applying part (1), we have a measurable set $E\subset \T$ with $m(E)=m(\T)=1$ such that
\beq\label{PMinv4}
m(\sigma_1((S,\zeta)))=m(S)
\eneq
for all Borel subsets $S\subset \T$ and $\zeta\in E.$
Thus, if $\zeta\in E\cap (\T\setminus F),$ by (\ref{PMinv4}), it is well known that
either
\beq\label{MINVt1}
\sigma_1(\xi, \zeta)=\xi g_1(\zeta)
\eneq
for all $\xi\in \T$ and for some $g_1(\zeta)\in \T$ or
\beq\label{MINVt2}
\sigma_1(\xi, \zeta)=\bar {\xi} g_1(\zeta)
\eneq
for all $\xi\in \T$ and for some $g_1(\zeta)\in \T.$

{\rm (3)} This part follows immediately  from (2).
By considering the subset $\{(1,\zeta): \zeta\in \T\}$ and applying (\ref{MINV3}) and
(\ref{MINV3+}),
we conclude that $\sigma_1(1, \zeta)=g_1(\zeta)$ is a continuous function.
Then,
by continuity of $\sigma, $ (3) follows.

\end{proof}

% Let $F\subset \T^2$
%be the closed subset so that $\sigma$ is continuous on $T^2\setminus F$ and $m_2(F)=0.$
%For each $\zeta\in \T,$
%$$
%\{\xi\in \T: (\xi, \zeta)\in F\}
%$$
%is a closed subset.  There exists a measurable subset $E_1\subset \T$ with $m(E_1)=m(\T)$
%such that
%\beq\label{PMinv5}
%m(\{\xi\in \T: (\xi, \zeta)\in F\})=0
%\eneq
%for all $\zeta\in E_1.$
%Put $E_2=E\cap E_1$ and $F(\zeta)=\{\xi\in \T: (\xi, \zeta)\in F\}.$
%Then $m(E_2)=m(\T)=1$ and $m(F(\zeta))=0$ for all $\zeta\in E_2.$
%
%Fix $\zeta\in E_2.$
%Then $\sigma_1(-, \zeta)$ is continuous on $O(\zeta)=\T\setminus F(\zeta).$
%Write $O(\zeta)=\cup_n O_n(\zeta),$ where $\{O_1,O_2,...\}$ is a set of pairwise disjoint
%open arks. By (\ref{PMinv4}), on each $O_n(\zeta),$ one must have
%\beq\label{PMinv6}
%\sigma_1((\xi, \zeta))=\xi c_n(\zeta)\,\,\,{\rm or}\,\,\, \sigma_1((\xi,\zeta))=\overline{\xi}c_n(\zeta)
%\eneq
%for some constant $c_n(\zeta)\in \T$
%Define $g_1(\xi, \zeta)=c_n(\zeta)$ for $\xi\in O_n(\zeta)$ and $\zeta\in E_2.$
%One can extend $g_1$ to a Borel function on $\T^2$  with values in $\T.$

{\bf The proof of Theorem \ref{MT1+}}
\begin{proof}

Let $\sigma((\xi, \zeta))=({\bar \xi}, \zeta).$ Then $\sigma: \T^2\to \T^2$ is a homeomorphism and
$\sigma^{-1}=\sigma.$
One has
\beq\label{DM1T1}
\sigma^{-1}\circ \Phi_{\theta, 1,0}\sigma((\xi, \zeta))
&=&\sigma^{-1}((\overline{\xi}e^{i2\pi\theta}, \zeta\overline{\xi}))\\
&=&(\xi^{-i2\pi\theta}, \zeta\xi^{-1})
\eneq
for all $\xi, \zeta\in \T.$
It follows that $\Phi_{\theta,1,0}$ and $\Phi_{-\theta, -1,0}$ are conjugate.

We will show that if $\overline{\theta_1\pm \theta_2}=0$ in $\R/\Z,$ then (3) in
\ref{MT1+} holds.  For the convenience, we will say $\af$ and $\bt$ are approximately conjugate in
the sense (M2) if  (3) holds. In this part of the proof, we will identify
$\T$ with
$\R/\Z.$

Let $\theta\in (0,1)$ be an irrational number. Since $\Phi_{\theta,1,0}$ and $\Phi_{-\theta, -1,0}$
are conjugate as shown above, it suffices to show that $\af$ and $\bt$ are approximately conjugate in
the sense (M2) if $\af=\Phi_{\theta, d_1, f_1}$ and $\bt=\Phi_{\theta, d_2, f_2}.$

Let $\ep>0.$ Choose $n>0$ so that $1/n<\ep.$ Let $J_1,J_2,...,J_k$ be the open arks provided by Lemma
\ref{Rok2} with the integer $n$ (and
$\af(t)=t+\theta$ for $t\in \R/\Z$)   Put $F(\xi)=\xi^{d_1}\exp(i2\pi f_1(\xi))$ and
$G(\xi)=\xi^{d_2}\exp(i2\pi f_2(\xi))$ for
$\xi\in \T.$
Let $\omega$ be the function in Lemma \ref{Lconj} with $\af(t)=t+\theta$ and $I_i=J_i,$ $i=1,2,...,k.$
By extending $\omega$ continuously on the left-sided closed arks $J_i^{(s)}$ for each $s$ and $i,$
as the last part of Lemma \ref{Lconj}, one has
\beq\label{PM1T1}
|[F(x)+\omega(x)]-[\omega(x+\theta)+G(x)]|\le \ep
\eneq
for all $x\in \R/\Z.$
Define
$$
\sigma((x,t))=(x,t+\omega(x))\,\,\,{\rm for}\,\,\, t, x\in \R/\Z.
$$
Therefore,
$$
m_2(\sigma(S))=m_2(S)
$$
for all Borel set $S\subset \T^2.$

Moreover, $\sigma$ is continuous except  finitely many circles (with the form
of $B\times \T,$ where $B$ is a finite subset of $\T$).

Now,
by (\ref{PM1T1}),
\beq\label{PM1T2}
\di(\af\circ \sigma((x,t)),\sigma\circ \bt((x,t)))\le |[F(x)+\omega(x)]-[\omega(x+\theta)+G(x)]|<\ep
\eneq
for all $x,t\in \R/\Z.$

For the converse, let $\theta_1$ and $\theta_2$ be two irrational numbers such that
$\overline{|\theta_1\pm \theta_2|}\not=0$ in $\R/\Z.$

Since we have shown that $\Phi_{\theta_1, d_1, f_1}$ and $\Phi_{\theta_1,1,0}$ are  approximately
conjugate in the sense (M2) above, and $\Phi_{\theta_2,d_2, f_2}$ and
$\Phi_{\theta_2,1,0}$ are approximately conjugate in the sense (M2) above, respectively,
it suffices to show that $\Phi_{\theta_1,1,0}$ and $\Phi_{\theta_2,1,0}$ are not approximately
conjugate in the sense (M2).

Put $\af=\Phi_{\theta_1,1,0}$ and $\beta=\Phi_{\theta_2,1,0}.$
Put
\beq\label{PMt0}
a=|e^{2\pi (\theta_1-\theta_2)}-1|>0 \andeqn b=|e^{2\pi (\theta_1+\theta_2)}-1|>0.
\eneq
Let $\ep>0$ such that $\ep<\min\{a/2,b/2\}.$ Suppose that there exists $\sigma: \T^2\to \T^2$ such
that
\beq\label{PMt1+}
\sup\{\di(\af\circ \sigma(x), \sigma\circ \bt(x)): x\in \T^2\}<\ep
\eneq
and $m_2(\sigma(S))=m_2(S)$ for all Borel  sets $S\subset \T^2.$

By  \ref{MINV}, there exists a Borel set $E\subset \T$ such that
$m(\T\setminus E)=0$ and  for each $\zeta\in E,$
\beq\label{PMt2}
\sigma_1((\xi,\zeta))=\xi g_1(\zeta), \,\,\,{\rm or}\,\,\,
\sigma_1((\xi,\zeta))={\bar \xi}g_1(\zeta),
\eneq
for all $\xi\in \T.$
Thus, for $\zeta\in E,$ by (\ref{PMt1+}), we have
\beq\label{PMt3}
|\xi g_1(\zeta)e^{i2\pi\theta_1}
-\xi g_1(\zeta\xi)e^{i2\pi\theta_2}|<\ep,\,\,\,{\rm or}\,\,\,
\eneq
\beq\label{PMt3+}
|\overline{\xi}g_1(\zeta)e^{i2\pi\theta_1} -\overline{\xi} g_1(\zeta\xi)e^{-i2\pi\theta_2}|<\ep,
\,\,\,{\rm or}\,\,\,
\eneq
\beq\label{PMt4}
|\xi g_1(\zeta)e^{i2\pi\theta_1}
-{\bar \xi} g_1(\zeta\xi)e^{-i2\pi\theta_2}|<\ep\,\,\,{\rm or}\,\,\,
\eneq
\beq\label{PMt4+}
|\overline{\xi}g_1(\zeta)e^{i2\pi\theta_1}
-\xi g_1(\zeta\xi)e^{i2\pi\theta_2}|<\ep
\eneq
for all $\xi\in \T.$

Choose $\xi=1,$ for all $\zeta\in E,$
one computes that either
\beq\label{PMT3+}
a=|g_1(\zeta)e^{2i\pi\theta_1}-e^{2i\pi\theta_2}g_1(\zeta)|<\ep\,\,\,{\rm or}\,\,\,
\eneq
\beq\label{PMT4+}
b=|g_1(\zeta)e^{2i\pi\theta_1}-e^{-2i\pi\theta_2}g_1(\zeta)|<\ep
\eneq
By (\ref{PMt0}), neither is possible.

\end{proof}

{\bf The Proof of Theorem \ref{MT1}}

\begin{proof}

(1) $\Rightarrow$ (2): We will modify the relevant part of the proof of Theorem \ref{MT1+}.
 Let
$\ep>0.$ By the proof of Theorem \ref{MT1+}, there exists a finite subset $B\subset \T$ and
a function $\omega: \T\to \T$ which is continuous on $\T\setminus B$ such that
\beq\label{PM1t+1}
\di(\sigma\circ \af((x,t)), \bt\circ\sigma((x,t)))<\ep
\eneq
for all $x,t \in \R/\Z,$ where $\sigma((x,t))=(x,t+\omega(x))$ for all $x, t\in \R/\Z.$ There is an
open subset $G\subset \T$ such that $B\subset G$ and
\beq\label{PM1t+2}
m(G)<\ep
\eneq
There is a continuous function $\omega_0$ from $\T$ to $\R/\Z$ such that
\beq\label{PM1t+3}
\omega_0(x)=\omega(x) \,\,\, {\rm for}\,\,\, x\in \T\setminus G.
\eneq
Now define $\sigma_0((x, t))=(x,  t+\omega_0(x))$ for $x, t\in \R/\Z.$ Then
\beq\label{PM1t+4}
\{(x,t): \di(\sigma_0\circ \af((x,t)), \bt\circ\sigma_0((x,t)))\ge \ep\}\subset \T\times G.
\eneq
It follows that
\beq\label{PM1t+5}
m_2(\{(x,t): \di(\sigma_0\circ \af((x,t)),\bt\circ\sigma_0((x,t)))\ge \ep\})<\ep.
\eneq
This proves (1) $\Rightarrow$ (2).

To see (2) $\Rightarrow$ (1), let $\theta_1$ and $\theta_2$ be two irrational numbers such that
$\overline{|\theta_1\pm \theta_2|}\not=0$ in $\R/\Z.$

Since we have shown that $\Phi_{\theta_1, d_1, f_1}$ and $\Phi_{\theta_1,1,0}$ are  approximately
conjugate in the sense of (M1), and $\Phi_{\theta_2,d_2, f_2}$ and
$\Phi_{\theta_2,1,0}$ are approximately conjugate in the sense  (M1), respectively,
it suffices to show that $\Phi_{\theta_1,1,0}$ and $\Phi_{\theta_2,1,0}$ are not approximately
conjugate  in the sense (M1).

Put $\af=\Phi_{\theta_1,1,0}$ and $\beta=\Phi_{\theta_2,1,0}.$ Suppose that $\sigma:\T^2\to \T^2$ is a
homeomorphism such that
\beq\label{PMt1}
m_2(\sigma(S))=m_2(S).
\eneq
Write $\sigma((\xi,\zeta))=(\sigma_1(\xi, \zeta), \sigma_2(\xi, \zeta)).$ It follows from (3) of
\ref{MINV} that
%By (\ref{PMt1}), for almost all $\zeta,$ $\sigma_1(\xi,\zeta)$ must be
%the form
\beq\label{PMtt2}
{\rm either}\,\,\,\sigma_1(\xi,\zeta)=\xi g_1(\zeta), \,\,\,{\rm or}\,\,\,
\sigma_1(\xi,\zeta)={\bar \xi}g_1(\zeta)
\eneq
for all $\xi, \zeta\in \T,$
where $g_1: \T\to \T$ is a continuous map.
Put
\beq\label{MPt+1}
a=|e^{2\pi (\theta_1-\theta_2)}-1|>0 \andeqn b=|e^{2i\pi (\theta_1+\theta_2)}-1|>0.
\eneq
Let $0<\ep<\min\{1/4,a/4, b/4\}.$ Let $z\in C(\T)$ be defined by $z((\xi, \zeta))=\xi$ for $\xi\in
\T.$ If $\af$ and $\bt$ are approximately conjugate  in the sense (M1), then there exists $\sigma$
described above such that
\beq\label{PMt3-}
\int_\T\int_\T |z(\af\circ \sigma(\xi, \zeta)-z(\sigma\circ\bt(\xi,\zeta))|d\xi d\zeta<\ep^4/4.
\eneq
One then computes that
\beq\label{PMtt3}
\int_{\T}|z(\af\circ \sigma(\xi, \zeta)-z(\sigma\circ\bt(\xi,\zeta))|d\xi<\ep^2
\eneq
for all $\zeta\in E,$ where $E$ is a measurable set such that
$m(E)>1-\ep.$

We now assume that $\sigma_1(\xi, \zeta)=\xi g_1(\zeta)$  for all $\xi, \zeta\in \T.$

For $\zeta\in E,$
\beq\label{PMtt4}
\ep^2&>&\int_{\T}|\xi g_1(\zeta)e^{i2\pi \theta_1}-\xi e^{i2\pi \theta_2} g_1(\zeta\xi)|d\xi\\
&=&\int_{\T}|g_1(\zeta)e^{i2\pi (\theta_1-\theta_2)}-g_1(\zeta\xi)|d\xi.
\eneq
%Fix $\zeta\in E,$ by (\ref{PMt4}),
%\beq\label{PMt5}
%\int_{\T}|g_1(\zeta_0)e^{i2\pi (\theta_1-\theta_2)}-g_1(\zeta_0\xi)|d\xi<\ep.
%\eneq
 By considering the constant function $F_1(\xi)=g_1(\zeta)e^{i2\pi (\theta_1-\theta_2)}$
and the function $F_2(\xi)=g_1(\zeta\xi)$ and by translating by ${\bar \zeta},$ we obtain
\beq\label{PMt6}
\int_{\T}|g_1(\zeta)e^{i2\pi (\theta_1-\theta_2)}-g_1(\xi)|d\xi<\ep^2
\eneq
for all $\zeta\in E.$  Fix $\zeta_0\in E$ and let
$$
E_{\zeta_0}=\{\xi\in |g_1(\zeta_0)e^{i2\pi (\theta_1-\theta_2)}-g_1(\xi)|<\ep\}.
$$
Then by (\ref{PMt6}), we compute that $m(E_{\zeta_0})>1-\ep.$ Therefore $m(E_{\zeta_0}\cap E)>0.$
 If
$\zeta_1\in E_{\zeta_0}\cap E,$ by (\ref{PMt6}), we have  that
\beq\label{PMt7}
&&\hspace{-0.8in}\ep^2 > \int_{\T}|g_1(\zeta_1)e^{i2\pi (\theta_1-\theta_2)}-g_1(\xi)|d\xi\\
&&\hspace{-0.7in} \ge  |g_1(\zeta_0)e^{i2\pi(\theta_1-\theta_2)}-g_1(\zeta_0)e^{i4\pi(\theta_1-\theta_2)}|\\
&&\hspace{-0.4in}-|g_1(\zeta_0)e^{i4\pi(\theta_1-\theta_2)}-g_1(\zeta_1)e^{i2\pi(\theta_1-\theta_2})|
-
\int_{\T}|g_1(\zeta_0)e^{i2\pi (\theta_1-\theta_2)}-g_1(\xi)|d\xi\\
&&\hspace{-0.7in}> |e^{i2\pi(\theta_1-\theta_2)}-1|-\ep-\ep^2>a/2.
\eneq
By the choice of $\ep$,  this is impossible.

Now we assume that $\sigma_1(\xi, \zeta)=\overline{\xi} g_1(\zeta)$ for all
$\xi, \zeta\in \T.$
As above, one has
\beq\label{PMt9}
\ep^4&>&\int_{\T}|{\overline{\xi} }g_1(\zeta)e^{i2\pi \theta_1}-\overline{\xi }e^{-i2\pi \theta_2}g_1(\zeta\xi)|d\xi\\
&=&\int_{\T}|g_1(\zeta)e^{i2\pi (\theta_1+\theta_2)}-g_1(\zeta\xi)|d\xi.
\eneq
The same argument used above leads us to
\beq\label{PMt10}
\ep^2> |e^{i2\pi (\theta_1+\theta_2)}-1|-\ep-\ep^2>{b\over{2}}.
\eneq
This would violate the choice of $\ep.$

\end{proof}

\section{Approximate  $K$-Conjugacy}

%{\bf Proof of Theorem{MT1}} Let $u$ be a unitary in $A_\af$ such that
%$uj_\af(f)u^*=j_\af(f\circ \af)$ for $f\in C(\T^2)$ and $v\in A_\bt$ be unitary such that
%$vj_\bt(f)v^*=j_\bt(f\circ \bt)$ for $f\in C(\T^2).$

%Let $\Lambda_n,$ $\Gamma_n$ and $\Sigma_n$ be as in \ref{IIIT1}.
%Let $\phi: A_{\theta}\to A_\af$ be the \hm\, defined in ?.
%Let $p\in A_\theta$ be the projection such that $\tau(\phi(p))=\theta,$ where
%$\tau$ is the unique tracial state on $A_{\theta}.$
%We may write $p=j(f_0) +j(f_1)u+u^*j(f_1),$ where $f_0, f_1\in C(\T)$ are non-negative functions, and
%where
%$j: C(\T)\to A_\theta$ is the embedding induced by rotation. Thus we obtain continuous functions $g_0,
%g_1\in T(\T^2)$ such that
%$g_0(x,t)=f_0(x)$ and $g_1(x,t)=f_1(x)$ so that
%$j_\af(g_0)=\phi(j(f_0)$ and $j_\af(g_1)=\phi(j(f_1)).$
%Then we have
%\beq\label{mt11}
%\Phi_n(\phi(p))=j_\bt(g_0\circ \Lambda_n)+j_\bt(g_1\circ \Lambda_n))v+v^*j_\bt(g_1\circ \Lambda_n).
%\eneq
%\andeqn
%Let $q=j_\bt(g_0)+j_\bt(g_1)v+v^*j_\bt(g_1).$ We have $\tau(q)=\theta.$ Moreover, by (\ref{t5}),
%\beq\label{mt12}
%\lim_{n\to\infty}|\Phi_n(\phi(p))-q\|=0.
%\eneq
%It follows that $\tau([\Phi_n(\phi(p))])=\theta$ for all large $n.$ This also applies to
%$\{\Psi_n\}.$ It follows that $\af$ and $\bt$ are approximately $\tau$-conjugate.

\vspace{0.2in}

.

\begin{Lem}\label{Ldense}
Let $\theta\in [0,1]$ be an irrational number and let
$$
V=\{a\sin kt+b\cos mt: a, b\in {\mathbb R}, k,m\in {\mathbb Z}, t\in [0,2\pi]\}.
$$
Then, for every $f\in V, $ there exists $g\in V$ such that
\beq\label{Ld1}
f(t)=g(t)-g(t+\theta)\,\,\rforal t\in [0,2\pi].
\eneq

\end{Lem}

\begin{proof}
Put
$$
V_0=\{ f(t)=g(t)-g(t+\theta):g\in V\}.
$$
It is clear that $V_0$ is a (real) vector space.

One has two elementary inequalities:  for any integer $n>1,$
\beq\label{Trig1}
|\sum_{k=1}^{n}\sin k\theta|&=&{|\cos {\theta\over{2}}-\cos (n+{1\over{2}})\theta|
\over{2\sin {\theta\over{2}}}}\\
&\le & {1\over{\sin {\theta\over{2}}}}
\eneq
and similarly
\beq\label{Trig2}
|\sum_{k=1}^n \cos k\theta| \le {1\over{\sin {\theta\over{2}}}}.
\eneq

Now, for any integer $m\in\Z,$
\beq\label{Trig3}
|\sum_{k=0}^n\sin (mt+k\theta)|&=&
|\cos mt (\sum_{k=0}^n \sin k\theta)+\sin mt (\sum_{k=0}^n \cos k\theta)|\\
&\le & |\sum_{k=0}^n \sin k\theta|+|\sum_{k=0}^n \cos k\theta)|\\
&\le & 1+{2\over{\sin {\theta\over{2}}}}
\eneq
for all $t\in \R.$
Now, since $t\mapsto t+\theta$ ($t\in \R/\Z$) is a minimal homeomorphism on $\T,$
by a lemma of Furstenberg (Lemma 5.2 of \cite{F}), there is $g\in C(\T)$ (real)
such that $\sin mt= g(t)-g(t+\theta)$ (for $t\in [0,2\pi]$).
It follows that $\sin mt\in V_0.$ Similarly, $\cos mt \in V_0.$
Since $V_0$ is a real vector space, $V\subset V_0.$
\end{proof}

The above  lemma can be proved directly by some trigonometric identities and
the function $g$ in the proof  may be chosen to be in $V.$

%Note that
%\beq\nonumber
%\sin kt-\sin(kt+k\theta)&=&(1-\cos k\theta)\sin kt+(-\sin k\theta)\cos kt\\
%&=& ({1-\cos k\theta\over{\cos \af}})\sin k(t+\af),
%\eneq
%where $\tan \af={-\sin k\theta\over{1-\cos k\theta}}.$ It follows that $({1-\cos k\theta\over{\cos
%\af}})\sin k(t+\af)\in V_0.$ Since $V_0$ is a vector space, $\sin k(t+\af)\in V_0.$ Similarly, $\cos
%k(t+\af)\in V_0.$ Since
%$$
%\sin kt=\cos k(-\af)\sin k(t+\af)+\sin k(-\af)\cos k(t+\af)
%$$
%and $V_0$ is a vector space,
%$$
%\sin kt\in V_0
%$$
%for all integer $k\in \Z.$ Similarly,
%$$
%\cos kt\in V_0.
%$$
%The  lemma follows.

\begin{Lem}\label{Ldense2}
Let $\theta\in [0,1]$ be an irrational number and let $f\in C(\T)$ be a real function. Then, for any
$\ep>0, $ there exists a continuous map $g:\T\to \T$ such that
$$
|\exp(i2\pi f(\xi))g(\xi){\overline{g(\xi\cdot e^{i 2\pi\theta})}}-1|<\ep
$$
for all $\xi\in \T.$

Moreover, if $d\not=0$ is an integer, $g$ may be chosen to have the form
$$g(\xi)=\xi^{kd}\exp (i2\pi g_0(\xi)),$$ where $k\in \Z$ and $g_0\in C(\T)$ is a real function.
% and where
%$V$ is defined in \ref{Ldense}.

\end{Lem}

\begin{proof}
Note that $\Z+\Z(\theta)$ is dense in $\R.$
Thus $[\Z+\Z(\theta)]/\Z$ is dense in $\R/\Z.$
Therefore, for any $a\in \R,$ there exists an integer $k\in \Z$ such that
\beq\label{Ld2}
|e^{i2\pi a}e^{-i2\pi k\theta}-1|<\ep
\eneq
Hence
\beq\label{Ld3}
|e^{i2\pi a}e^{i2\pi kt}e^{-i2\pi k(t+\theta)}-1|&=&|e^{i2\pi a}e^{-i2\pi k\theta}-1|<\ep
\eneq
for all $t\in [0,1].$

Let $V$ be as in the proof of \ref{Ldense} and let $f_0\in V.$
Applying \ref{Ldense} and choose a real $g_0\in C(\T)$ such that
\beq\label{Ld4}
f_0(t)=g_0(t+\theta)-g_0(t)
\eneq
for all $t\in \R/\Z.$
Let $f=a+f_0$ and $g(t)=\exp(i2\pi(kt+g_0(t)))$ for $t\in \R/\Z.$
Then
\beq\label{Ld5}
&&\hspace{-0.4in}|\exp(i2\pi f(t))g(t){\overline{g(t+\theta)}}-1|\\
&&\hspace{-0.4in}=|\exp\{i2\pi f(t)\}\exp\{i2\pi(kt+g_0(t))\}\exp\{-i2\pi [k(t+\theta)+g_0(t+\theta)]\}-1|\\
&&\hspace{-0.4in}=|\exp\{i2\pi (a+f_0)\}\exp\{i2\pi(-k\theta +g_0(t)-g_0(t+\theta))\}-1|\\
&&\hspace{-0.4in}=|\exp\{i2\pi(a-k\theta)\}-1|<\ep
\eneq
for all $t\in \R/\Z.$

By the Stone-Weierstrass theorem, the set of real trigonometric polynomials is dense in the real part
of $C(\T).$ Thus the first part of the lemma follows.

To see the last part of the lemma, we only need to note that $d\theta$ is also an irrational number
and $\Z d\theta$ is dense in $\R/\Z.$
\end{proof}

{\bf Proof of \ref{MT2}}

\begin{proof}

That (2) $\Leftrightarrow$ (4) follows the computation in Example 4.9 of \cite{Ph} (see \ref{RCK}) and
(1) $\Leftrightarrow$ (4) follows from the classification theorem in \cite{Lnduke} as mentioned at the
end of \ref{Rdkc}. It is also clear that (3)
$\Rightarrow$ (5).

It remains to show (2) $\Rightarrow$ (3) and (5) $\Rightarrow$ (2).

We will first show (2) $\Rightarrow$ (3).

 As in the proof of \ref{MT1}, $\Phi_{\theta, d, 0}$ and $\Phi_{-\theta, -d, 0}$ are
conjugate. Define
$\sigma((\xi, \zeta))=(\xi, \overline{\zeta}).$ Then $\sigma: \T^2\to \T^2$ is a
homeomorphism and $\sigma^{-1}=\sigma.$  One check that
\beq\label{M2T1}
\sigma^{-1}\circ \Phi_{\theta, d, 0}\circ \sigma((\xi, \zeta))&=&
\sigma^{-1}((\xi e^{2i\pi\theta}, \overline{\zeta} \xi^{d}))\\
&=&(\xi e^{2i\pi\theta},\zeta \xi^{-d})=\Phi_{\theta, -d,0}((\xi, \zeta))
\eneq
for all $\xi, \zeta\in \T.$ Therefore
$\Phi_{\theta, d,0}$ and $\Phi_{\theta, -d,0}$ are conjugate.
So $\Phi_{-\theta, -d, 0}$ and $\Phi_{-\theta, d,0}$ are conjugate. Combining with the fact that
$\Phi_{\theta, d, 0}$ and $\Phi_{-\theta, -d, 0}$ are conjugate that mentioned above,
we conclude that
$\Phi_{\theta, d,0}$ and $\Phi_{-\theta, d,0}$ are conjugate.
Thus, to complete the proof, it suffices to show that $\af=\Phi_{\theta, d, f_1}$ and
$\bt=\Phi_{\theta, d, f_2}$ are approximately $K$-conjugate for any real continuous Lipschitz
functions $f_1$ and $f_2.$

It follows from Theorem 4.6 of \cite{LP} that both $A_\af$ and $A_\bt$ have tracial rank zero. By the
$K$-theory computation in \ref{RCK}, there is an order isomorphism
$$
\kappa: (K_0(A_\af), K_0(A_\af)_+, [1_{A_\af}], K_1(A_\af))
\to (K_0(A_\bt), K_0(A_\bt)_+, [1_{A_\bt}], K_1(A_\bt))
$$
such that $\kappa([u_\af])=[u_\bt].$ By the classification theorem (\cite{Lnduke}, there exists a
unital isomorphism $\phi: A_\af\to A_\bt$ such that
$$
[\phi]=\kappa.
$$

Let $f(\xi)=f_2(\xi)-f_1(\xi)$ for $\xi\in \T.$ Fix $\dt>0.$ By applying \ref{Ldense2}, we obtain
\beq\label{PMTT1}
g(\xi)=\xi^{kd} \exp (2i\pi g_0(\xi))
\eneq
for $\xi\in \T,$ where $g_0\in C(\T)$ is a real function such that
\beq\label{PMTT2}
|\exp (2i\pi f(\xi))g(\xi){\overline{g(\xi e^{2i\pi \theta})}}-1|<\dt
\eneq
for all $\xi\in \T.$ Define
\beq\label{PMTT3}
\sigma((\xi, \zeta))=(\xi, \zeta g(\xi))
\eneq
for all $(\xi, \zeta)\in \T^2.$ Then
\beq\label{PMTT3+}
&&\sigma\circ \af((\xi, \zeta))=(\xi e^{2i\pi\theta}, \zeta \xi^d e^{2i\pi f_1(\xi)}g(\xi
e^{2i\pi\theta}))\andeqn\\
&&\bt\circ \sigma((\xi, \zeta))=(\xi e^{2i\pi\theta}, \zeta g(\xi) \xi^d e^{2i\pi f_2(\xi)})
\eneq
for all $\xi, \zeta\in \T.$  Using (\ref{PMTT2}), we estimate that
\beq\label{PMTT4}
\di(\sigma\circ \af((\xi, \zeta)),\bt\circ\sigma((\xi, \zeta)))<\dt
\eneq
for all $(\xi, \zeta)\in \T^2.$

Note that $\sigma$ is homotopic to $\Phi_{\theta, kd, g_0}.$ As the computation in Example
4.9 of \cite{Ph} (see \ref{RCK}), we have
$\sigma_{*0}={\id}_{*0}$ on $K_0(C(\T^2))$ and
$\sigma_{*1}$ on $K_1(C(\T^2))\cong \Z\oplus \Z$ is represented by the matrix
$$
\begin{pmatrix}1 & kd\cr
                        0  & 1\cr
                        \end{pmatrix}.
                        $$

It induces the identity map from $\Z/d\Z\oplus \Z$ onto $\Z/d\Z\oplus \Z.$ Therefore
\beq\label{PMTT5}
h_{*i}=(\phi\circ j_\af)_{*i}=(j_\bt)_{*i}, \,\,\,i=0,1,
\eneq
where $h: C(\T^2)\to A_\bt$ is defined by $h(f)=\phi\circ j_\af(f\circ \sigma)$ for
$f\in C(\T^2).$ Since $K_i(C(\T^2))$ is free, we, in fact, have that
\beq\label{PMTT6}
[h]=[\phi\circ j_\af]=[j_\bt]\,\,\,{\rm in}\,\,\, KL(C(\T^2), A_\bt).
\eneq
 We also note that
 \beq\label{PMTT7}
 \tau\circ h(f)=\tau\circ \phi\circ j_\af(f)=\tau\circ j_\bt(f)
 \eneq
for all $f\in C(\T^2),$ where $\tau$ is the unique tracial state on $A_\bt.$

Let
$${\cal F}_1={\cal F}\cup
\{f\circ\sigma^{-1}: f\in {\cal F}\}\cup\{f(\sigma\circ \af\sigma^{-1}): f\in {\cal F}\}.
$$
By (\ref{PMTT6})  and (\ref{PMTT7}), and by Theorem 3.4 of
\cite{LnKC}, there exists a unitary $W\in A_\bt$
such that
\beq\label{M2T101}
Wj_\bt(f)W^*\approx_{\ep/3} \phi\circ j_\af(f\circ \sigma)\,\,\, {\rm on}\,\,\, {\cal F}_1.
\eneq
In particular, if $f\in {\cal F},$
\beq\label{M2T102}
W^*\phi\circ j_\af(f) W\approx_{\ep/3} j_\bt(f\circ \sigma^{-1})\andeqn
\eneq
\beq\label{M2T103}
W^*\phi\circ j_\af(f(\sigma\circ \af) W\approx_{\ep/3} j_\bt(f(\sigma\circ\af\circ \sigma^{-1})).
\eneq
Therefore
\beq\label{M2T104}
W^*\phi(u_\af^*)Wj_\bt(f)W^*\phi(u_\af)W  &\approx_{\ep/3}& W^*\phi(u_\af^*)\phi\circ
j_\af(f\circ\sigma)\phi(u_\af)W\\
&=&W^*\phi\circ j_\af(f\circ\sigma\circ\af)W\\
&\approx_{\ep/3}&j_\bt(f\circ \sigma\circ\af\circ\sigma^{-1})
\eneq
for all $f\in {\cal F}.$ It follows that, with sufficiently small $\dt,$
\beq\label{M2T105}
{\rm ad}\, (W^*\phi(u_\af)W)\circ (j_\bt(f))\approx_{\ep} j_\bt(f\circ\bt)
\eneq
for all $f\in {\cal F}.$ Put $z=u_\bt^*( W^*\phi(u_\af)W).$ Then, since $[\phi(u_\af)]=[u_\bt],$
$z\in U_0(A_\bt)$ and
\beq\label{M2T106}
j_\bt(f)z\approx_{\ep} zj_\bt(f)
\eneq
for all $f\in {\cal F}.$ This shows that $\af$ and $\bt$ are approximately $K$-conjugate.

Now we consider (5) $\Rightarrow$ (2). Suppose that $\af=\Phi_{\theta_1, d_1, f_1}$ and
$\bt=\Phi_{\theta_2,d_2, f_2}$ are two Furstenberg transformations.
Suppose that there exist sequences of homeomorphisms $\{\sigma_n\}$ and $\{\gamma_n\}$ on $\T^2$ such
that
\beq\label{PMTTT1}
&&\hspace{-0.6in}\lim_{n\to\infty}\sup\{ \di(\sigma_n\circ \af\circ \sigma_n^{-1}((\xi,\zeta)),
\bt((\xi,
\zeta))): (\xi, \zeta)\in
\T^2\}=0\\
&&\hspace{-0.6in}\lim_{n\to\infty}\sup\{\di(\gamma_n\circ \bt\circ\gamma_n^{-1}((\xi,
\zeta)),\af((\xi,
\zeta))):
(\xi, \zeta)\in \T^2\}=0,\\
&&\hspace{-0.4in}m_2(\sigma_n(S))=m_2(S)\andeqn m_2(\gamma_n(S))=m_2(S)
\eneq
for all Borel sets $S\subset \T^2.$ It follows from \ref{MT1} that
$\overline{\theta_1\pm \theta_2}=0$ in $\R/\Z.$

Write
$$
\sigma_n((\xi, \zeta))=(G_1^{(n)}((\xi, \zeta)), G_2^{(n)}((\xi, \zeta))
$$
for all $(\xi, \zeta)\in \T^2.$
 It follows from \ref{MINV} that there are continuous maps $g_1, g_2: \T\to \T$ such that
\beq\label{PMTTT2}
G_1^{(n)}((\xi, \zeta))={\tilde \xi} g_1^{(n)}(\zeta)\andeqn G_2^{(n)}((\xi,\zeta)) ={\tilde \zeta}
g_2^{(n)}(\xi)
\eneq
for all $\xi, \zeta\in \T,$ where ${\tilde \xi}=\xi$ for all $\xi\in \T$ or
${\tilde \xi}={\overline{\xi}}$ for all $\xi\in T,$
${\tilde \zeta}=\zeta$ for all $\zeta\in \T$ or
${\tilde \zeta}={\overline{\zeta}}$ for all $\zeta\in \T.$
We have
\beq\label{PMTTT3}
&&\hspace{-0.4in}\sigma_n\circ \af((\xi, \zeta))=({\tilde \xi}e^{\pm
2i\pi\theta_1}g_1^{(n)}(\zeta\xi^{d_1}e^{2i\pi f_1(\xi)}),
{\tilde \zeta}\xi^{\pm d_1}e^{\pm i\pi f_1(\xi)} g_2^{(n)}(\xi e^{2i\pi\theta_1}))\\
&&\andeqn\nonumber\\
 &&\hspace{-0.4in}\bt\circ\sigma_n((\xi, \zeta)) =({\tilde
\xi}g_1^{(n)}(\zeta)e^{2i\pi\theta_2},{\tilde
\zeta}g_2^{(n)}(\xi)
\xi^{\pm d_2}g_1^{(n)}(\zeta)^{d_2}
e^{2i\pi f_2({\tilde \xi}g_1^{(n)}(\zeta))})
\eneq
for all $\xi,\zeta\in \T.$

We compute that, for all sufficiently large $n,$
\beq\label{PMTTT5}
|\xi^{\pm d_1\pm d_2} e^{2i \pi(f_2({\tilde \xi}g_1^{(n)}(\zeta))\pm f_1(\xi))}
g_2^{(n)}(\xi){\overline{g_2^{(n)}(\xi e^{2i\pi\theta})}} g_1^{(n)}(\zeta)^{d_2}-1|<1/2
\eneq
for all $\xi, \zeta\in \T.$ Fix $\zeta\in \T$ and let $\xi$ vary in $\T.$ It follows that, for fixed
$\zeta\in \T,$
$$
\xi^{\pm d_1\pm d_2} e^{2i \pi(f_2({\tilde \xi}g_1^{(n)}(\zeta))\pm f_1(\xi))}
g_2^{(n)}(\xi){\overline{g_2^{(n)}(\xi e^{2i\pi\theta})}} g_1^{(n)}(\zeta)^{d_2}
$$
is homotopically trivial as a unitary in $C(\T).$ Since
$g_2^{(n)}(\xi){\overline{g_2^{(n)}(\xi e^{2i\pi\theta})}}$ is homotopically trivial,
(\ref{PMTTT5}) implies that
%$g(\zeta)^{d_2}$
$\xi^{\pm d_1\pm d_2}$ is homotopically trivial. However, that can only happen
when $|d_1|=|d_2|.$

%It follows from Theorem 3.2 of \cite{LM1} that there exists a sequence of \morp s
%$\{\psi_n\}$ from $A_\bt$ to $A_\af$ such that

%(i) $\lim_{n\to\infty}\|\psi_n(ab)-\psi_n(a)\psi_n(b)\|=0$ for all $a, b\in A_\bt;$

%(ii) $\lim_{n\to\infty}\|\psi_n(u_\bt)-u_\af\|=0$ and

%(iii) $\lim_{n\to\infty}\|\psi_n\circ j_\bf(f)-j_\af(f\circ \sigma_n)\|=0$ for all $f\in C(\T^2).$

%Denote by $z$ the function on $C(\T^2)$ defined by $z((\xi, \zeta))=\xi.$ Then

%It follows from 3.3 of \cite{LM1} that for all large $n,$ there are \hm s from $\Z/d_2\Z$ maps to
%$\Z/d_1\Z$ and it maps the generator $[u_\bf]$ of $\Z/d_2\Z

\end{proof}

\begin{Cor}\label{C1}
Let
$$
V_1=\{m_1\theta+m_2 +a\cos m_2 t +b \sin m_4 t: m_1, m_2, m_3, m_4\in \Z,
a, b\in \R, t\in [0,2\pi]\}.
$$
Let $\af=\Phi_{\theta, d, f_1}$ and $\bt=\Phi_{\theta, d, f_2}$ such that
$f_1-f_2\in V_1.$
Then $\af$ and $\bt$ are conjugate.
Note that $V_1$ is dense in the real part of $C(\T).$
\end{Cor}

\end{document}